\input amstex
\input xy
\xyoption{all}
\documentstyle{amsppt}
\document
\magnification=1200
\NoBlackBoxes
\nologo
\hoffset1.5cm
\voffset2cm
\vsize16cm

\def\N{{\Bbb N}}
\def\S{{\Bbb S}}
\def\Q{\bold{Q}}

\def\Z{{\Bbb Z}}

\centerline{\bf SYMMETRIES OF GENUS ZERO MODULAR OPERAD}

\bigskip

\centerline{\bf No\'emie C. Combe and Yuri I. Manin}

\medskip

{\it ABSTRACT.} In this article combining survey and certain research results,
we introduce a categorical framework
for description of symmetries of genus zero modular operad. This
description merges the techniques of recent  ``persistence homology'' studies
and the classical formalism of groupoids. We show that the contravariant
 ``poset in groupoids'' embodying these symmetries, provides a new avatar
of profinite Grothendieck--Teichm\"uller group acting  upon this operad
but seemingly not related with representations of the Galois group
of all algebraic numbers.

\bigskip

\centerline{\bf 1.  Introduction and summary}

\medskip

Recently considerable attention was attracted by the studies of interaction
of modular operad (playing the central role in quantum cohomology constructions)
with the celebrated Grothendieck--Teichm\"uller group ([HaLoSch], [Ho], [BrHoRo], [Fr], [Lo]).
\smallskip

In our  paper [CoMa] motivated by these studies,
we started looking more
closely upon symmetries of the genus zero components of this operad.
We focused our attention there upon a  very particular symmetry of order two and
the respective quotient operad.

\smallskip

Here we take into account the whole family of automorphisms of $\overline{M}_{0,n}$, 
calculated in [BrMe]
and combine them into two  ``locally finite'' combinatorial objects (posets in groupoids)
reflecting their compatibilities
with operadic compositions (here results of [BrMe] are also crucial; their proofs are based
upon Kapranov's constructions from [Ka]).

\smallskip

This leads us to two versions of the notion of  ``full automorphism group'' of the genus zero
operad and opens  a question about possibility to transfer in this combinatorial
context also an action of the absolute Galois group ([I], [Sch]).
\smallskip

The point is that each study achieving this goal
which we are aware of,
passes through a homotopic environment,
and we wonder whether it is necessary. The dendroid formalism
used in [CiMoe] works with simplicial sets (or simplicial objects of other
categories), and our formalism here may be considered as a simplified version of
 their constructions. We avoid in this way simplicial constructions and 
hopefully, can thus bypass  the $\infty$--versions 
of the classical Grothendieck--Teichm\"uller correspondence. Perhaps,
a price we have to pay for this is an interaction with the Galois group which
gets ``lost in translation''.

\smallskip

Sections 2--5 of this article are short introductions to the combinatorics
of components of the genus zero modular operad and of the relevant operadic
compositions. Section 6 focusses on the main content of this note: a brief explanation of the structures
of two {\it"graphic skeleta''} of the symmetries of this operad, one of 
which embodies covariant symmetries and another contravariant ones.

\bigskip

\centerline{\bf 2. Cofinite combinatorics and operadic symmetries}

\medskip

Denote by $\N$ the set of natural numbers $\{1,2,3, \dots \}$. Call a subset
$S\subset \N$ {\it cofinite} one (cf for brevity), if its complement is finite.
Call a map $f:\, \N \to \N$ {\it cf--map}, if it is identical on an appropriate
cf--set (that may depend on $f$).

\smallskip

{\bf Claim.} {\it (a) Composition of two cf--maps is a cf--map.

(b) Bijective cf--maps form a group wrt composition.}

\smallskip
Denote the latter group $\S_{\infty}$. Clearly, it coincides
with the union $\cup_{n=1}^{\infty} \S_n$, in which
every symmetric group $\S_n$ is the full group of bijective
cf--maps identical outside $\{1,\dots , n\}$.

\medskip

Below, we work over a field of characteristic zero.

\medskip

{\bf 2.1. Theorem.}  {\it (a) $\S_{\infty}$ acts upon the family of operadic components $\overline{M}_{0,n}$,
$n\ge 3$, in the following way.

\smallskip

Let  $C_{0,n}\to \overline{M}_{0,n}$ be the universal family of stable curves of genus zero
endowed with $n$ structure sections $s_i:\, \overline{M}_{0,n} \to C_{0,n}$. For each $g\in \S_n$, produce
from it another family $g^{-1}(C)_{0,n}$ by renumbering the sections: $s_i$ acquires
the new marking $s_{g(i)}$. By universality, we obtain for an appropriate
automorphism $g_M: \,\overline{M}_{0,n}  \to  \overline{M}_{0,n}$ that the same
renumbering can be obtained via $g_M^*$.
\smallskip

This map $\S_n \to Aut\,\overline{M}_{0,n}$ is surjective for all $n \ge 3$ and bijective for $n\ge 5$.
\smallskip

(b) This family of automorphisms $\{ g_M\}$ can be naturally and uniquely extended to the family
of all operadic compositions between the components.}

\medskip

{\bf Sketch of proof.} (a) For a proof of the first statement, see [BrMe], Theorem 4.3.

\smallskip

(b) Regarding the second statement, we start with recalling  that the structure of
genus zero modular operad is given by a family of composition morphisms
$$
\overline{M}_{0,n_1} \times \dots \times \overline{M}_{0,n_r} \to \overline{M}_{0,n_0}.
$$
On the level of geometric fibers of universal families of curves, each such
morphism corresponds to a controlled degeneration of the respective fiber
of $C_{0,n_0}$ with marked points by grafting to it fibers of $C_{0,n_1}, \dots , C_{0,n_r}$,
and subsequent renumbering of all marked points, that have not been grafted.

\smallskip

The first part of this procedure is pure geometric one. Therefore it only remains to check
the compatibility of these ``operadic renumberings'' with combinatorics
of bijective renumberings used in the definition of
$\S_{n_i} \to Aut\,\overline{M}_{0,{n_i}}$. 
\smallskip

The interested reader will find more details in Sections 4 and 6. 

\medskip

{\bf Warning on boundary cases.} For $n\ge 4$, $\overline{M}_{0,n}$ is a smooth
projective variety of dimension $n-3$; $\overline{M}_{0,3}$ is a point; and
for $n=0,1,2$ it is a stack, the classifying stack of the automorphism
groups of  $P^1$ with $0,1,2$ marked points respectively.

\smallskip

Studying componentwise symmetries of the modular
operad we will usually omit an explicit discussion of special cases involving small values
of $n$.
\bigskip

\centerline{\bf 3. Stable curves of genus zero and their graphs}

\medskip

Below we will use essentially the language of genus zero quantum cohomology
as it was introduced in [KoMa] and developed in various directions in many
other papers including combinatorial formalisms of [BoMa].

\smallskip

Recall that a {\it stable} genus zero algebraic curve with marked points
is a family $(C, x_1, \dots , x_n)$. Here $C$ is an algebraic curve, possibly reducible,
with smooth irreducible components of genus zero. If it has singular points,
they are double points in which two different components intersect transversally,
From the genus zero condition it follows that a sequence of irreducible components,
in which every two consecutive curves intersect, cannot form a cycle.
A finite set of smooth points $(x_1,\dots ,x_n)$ of $C$ are marked.
Stability means that the automorphism group of such a curve, fixing marked points,
is finite. Equivalently, each irreducible component intersecting 
with other components at $k\ge 0$ points, carries at least $3-k$ marked points.

\smallskip

If we forget precise positions of marked and singular points on their
components, the remaining combinatorial structure of $C$
is well encoded by the {\it dual graph $\tau$} of this curve.

\smallskip
Such a graph   consists of {\it vertices, flags and edges.} Vertices are points,
bijectively corresponding to irreducible components of $C$, edges are segments
of a line bijectively corresponding to singular points of $C$. 
Two vertices are connected by an edge iff the respective components of $C$
intersect at the respective singular point. Finally, smooth marked points
bijectively correspond to flags: such a flag is a ``half--segment'', whose one end is attached
to the vertex corresponding to the respective component, and other end
is free (we prefer not to include this free end into the set of vertices, unlike some
other authors working with similar formalisms). An edge then might be imagined
as the union of two flags whose free ends are now connected.
If our curve $C$ is irreducible, then its tree has only one vertex,
to which several flags are attached. Such trees are called {\it corollas.}

\smallskip

In order not to mix free flags with halves of edges, one may call
a free flag {\it a leaf}, or {\it a tail} as in [KoMa].

\smallskip

 For  wider contexts and more precise descriptions,
generalizable to curves of higher genera,
see [BoMa].

\smallskip

Among various types of labelings of sets, forming a graph, an important role
is played by {\it orientations of flags.} If two flags form an edge, their orientations
must agree. We often use oriented trees, in which one  free flag is chosen as a {\it root},
and all other leaves are oriented in such a way, that from each leaf there exists a unique
oriented path to the root.

\smallskip

If $\tau$ corresponds to a curve $(C, x_1,\dots , x_n)$ on which marked points
are labeled, say, by $1,\dots ,n$, then it it is convenient to use the same labels to
mark tails  of $\tau$.
\bigskip

\centerline{\bf 4. Stratifications of $\overline{M}_{0,n}$}

\medskip

Let $\tau$ be the (labeled)  graph of a stable connected curve $(C,x_1,\dots ,x_n)$.
Denote by $M_{0,\tau} \subset \overline{M}_{0,n}$ the moduli submanifold (or generally, substack)
parametrizing  all curves having the same graph $\tau$.  Its closure will be denoted
$\overline{M}_{0,\tau}$.

\smallskip

{\bf Examples.} a) If $\tau$ is a corolla, $M_{0,\tau}$ is the maximal stratum,
having dimenstion $n-3 = dim\, \overline{M}_{0,n}$. 

\smallskip

b) If $\tau$ has two vertices, they must be connected by an edge. The  possible
distributions of leaves between vertices bijectively correspond to (unordered)
partitions of $\{x_1, ... ,x_n\}$ into two subsets of cardinality  $\ge 3$ each,
and in additional choices of one point in each part. Codimension
of $\overline{M}_{0,\tau}$ is 1.

\smallskip

c)  There are precisely $n+1$ strata that are {\it points} (codimension $n-3$).
The respective trees are all isomorphic if one forgets labeling of leaves:
they are represented by the sequence of vertices $\{v_1, \dots ,v_{n-1}\}$
such that $v_i$ and $v_{i+1}$ are connected by an edge for $i\le n-2$;
besides, there is one leaf at each $v_2,\dots , v_{n-3}$ and two leaves
at $v_1$ and $v_{n-1}$ each. All in all, we get $n$ leaves, and $\S_n$
acts simply transitive upon their ordering in the trees.

\smallskip

 In order to state a general result in convenient form, we will remind
 the basic definitions of morphisms of stable trees, as they were stated
 in [KoMa], Sec. 6.6--6.8, but now in the language of their
 geometric realizations as above.

  \smallskip
 
 A morphism $f:\,\tau \to \sigma$ is determined by its covariant surjective action
 upon vertices $f_v:\, V_{\tau}\to V_{\sigma}$ and contravariant injective actions
 upon tails and edges:
 $$
 f^t:\, T_{\sigma}\to T_{\tau}, \quad f^e:\, E_{\sigma}\to E_{\tau}.
 $$
 Geometrically,  $f$ contracts edges from $E_{\tau}\setminus f^e(E_{\sigma})$
 and tails from $T_{\tau}\setminus f^t(T_{\sigma})$,
 compatibly with its action upon vertices.
 
 \smallskip
 
 An important  operation on graphs corresponding to operadic composition was called {\it gluing}
 in [KoMa], Sec. 6.6.4. Starting with two pairs consisting of a tree and its tail
 $(\tau_i, t_i)$, $i=1,2$, we produce a new tree $(\tau_1,t_1) * (\tau_2, t_2)$
 by connecting $t_1$ and $t_2$ into one new edge.
 
 \smallskip
 
 Then we have the following functoriality property of this operation: for
 any two morphisms $f_i:\, \tau_i\to \sigma_i$ not contracting $t_i$ we
 have an obvious morphism
 $$
 f_1 * f_2:\ (\tau_1,t_1) * (\tau_2, t_2) \to (\sigma_1, (f_1^t)^{-1}(t_1)) * (\sigma_2, (f_2^t)^{-1}(t_2)) .
 $$
 The general result about geometry of stratification of $\overline{M}_{0,n}$ now can be stated
 as follows ([Ke], [KoMa]). Let $\overline{M}_{0,\tau}$ be the moduli space of stable curves
 of genus zero whose marked points are bijectively labeled by  tails of $\tau$.
 Let ${M}_{0,\tau}$ be the locally closed subspace of it where respective curves have combinatorial type 
 exactly $\tau$.
 
 \medskip
 
 {\bf 4.1. Theorem.} {\it There is a unique family of stratifications of $\overline{M}_{0,\tau}$
 by strata of the types ${M}_{0,\sigma}$ and  behaving in the functorial
 way wrt a generating class of tree morphisms.}  
 
 \bigskip
 
\centerline{\bf 5. Thin categories and groupoids}
 
 \medskip
 
  {\bf  5.1. Thin categories and groupoids.}  From now on, our graphs
may be infinite. To avoid set--theoretic complications, we will work in a fixed small universe:
see [KashSch], Sec. 1.1.

\smallskip

We start with introducing some basic notions of ``persistence formalism''. Here we rely
upon [BuSSc];  for more details and
basic sources and applications see Introduction and References in [MaMar].

\smallskip

A {\it proset} (preordered set) is a set $P$ endowed with a binary relation $\le$
which is reflexive and transitive.
\smallskip
A {\it poset} (partially ordered set) is a set $\Cal{S}$ endowed with a binary relation $\le$
which is reflexive, transitive and anti--symmetric. A proset is a poset iff $X\le Y$
and $Y\le X$ imply $X=Y$. It follows that each proset has a canonically defined
quotient which is a poset: identify in $P$ pair of objects $X, Y$ for which we
have simultaneously $X\le Y$ and $Y\le X$.

\smallskip

A poset $\Cal{S}$ defines an oriented graph (generally infinite)
whose vertices are (marked by) $\Cal{S}$, and edges connect pairs $X,Y\in \Cal{S}$ with $X<Y$
oriented from $X$ to $Y$.
\smallskip
A category $\Cal{C}$ is called {\it thin} if for any two objects
$X$, $Y$, the set  $\roman{Hom}_{\Cal{C}}(X,Y)$  consists of $\le 1$ element, and 
if both $\roman{Hom}_{\Cal{C}}(X,Y)$ and $\roman{Hom}_{\Cal{C}}(Y,X)$
are non--empty, then $X=Y$.

\smallskip

It follows that all
automorphisms of  $X$, $Y$ act upon  $\roman{Hom}_{\Cal{C}}(X,Y)$ as identity, so $\Cal{C}$ is
equivalent to a category for which $\roman{Hom}_{\Cal{C}}(X,X)= \{id_X\}$ fo any object $X$,
which we will temporarily assume.
\smallskip
   For such a category, $\roman{Ob}\  \Cal{C}$ has a canonical structure of  a poset:
 $X\le Y$ iff $\roman{Hom}(X,Y)$ is nonempty. Conversely, each poset
 defines in this way a thin category in which morphisms are equivalence classes  of 
 {\it oriented paths} from $X$ to $Y$. Hence, describing a thin category, one can restrict oneself to
 an explicit description of only some set of {\it generating morphisms} and keep in mind that
 each diagram in a thin category is automatically commutative.
 
 \smallskip
 
 In a sense, a complementary class of small categories is the class of {\it groupoids.}
 Recall that a  category $\Cal{G}$ is a groupoid, if all morphisms in it are isomorphisms
 ([KashSch], p. 13).
 
 \smallskip
 
 These two classes can be naturally merged. 
 
 \medskip
 {\bf 5.2. Definition.} {\it We will call a poset in groupoids 
 a category  $\Cal{P}\Cal{G}$, satisfying the following two conditions:
 
 \smallskip
 
 (a)  For any object $X$, the full subcategory of  $\Cal{P}\Cal{G}$
 consisting of all objects isomorphic to $X$ is a groupoid.
 
 \smallskip
 
 (b) If  $X$ and  $Y$ are not isomorphic and $\roman{Hom} (X,Y)$ is non--empty,
 then  $\roman{Hom} (X,Y)$ has a single orbit with respect to the
 precomposition by the automorphism group of $X$ and postcomposition 
 by the automorphism group of $Y$.}
 
 \smallskip
 
 There is a natural universal  functor from posets in groupoids to 
 thin categories, identical on objects and identifying all morphisms in each
 non--empty  $\roman{Hom} (X,Y)$.

  \medskip
  
  {\bf  5.3. Thin category of stable  trees.}  For two such trees $\pi, \tau$ we
  put $\pi < \tau$ if $\pi =(\sigma, t_1)* (\tau, t_2)$ for appropriate $\sigma$ and $t_1, t_2$,
  where we return to notations of Sec.  4. Passing
  to categories, we remind that this binary relation produces the generating morphisms of 
  the respective category of trees: the remaining ones are obtained by composing them.
  
  \smallskip
  
 Theorem 4.1 offers  another equivalent thin category whose objects are strata
$\Cal{M}_{0,\pi}$ of all components of genus zero modular operads.

\smallskip

There are  versions of these posets that are posets in groupoids.

\smallskip

Namely, consider stable trees whose tails form themselves (or  are bijectively labeled by)
finite sets $T$ of our small universe. Then  trees of given topological type with
fixed $T$ form a groupoid  whose morphisms can be identified
with bijections $T\to T$.
 \bigskip
 
\centerline{\bf 6. Graphic skeleta of the symmetries}

\smallskip

\centerline{\bf of genus zero modular operad}

\medskip

The widest categorical context in which we can imagine symmetries of
the genus zero operad would involve all forests $\tau$ in a small
universe and respective spaces $\overline{M}_{0,\tau}$ and $M_{0,\tau}$.

\smallskip

As was explained above, in our ``economy class version'', we will restrict ourselves to corollas,
that is finite sets $T$ marking points on stable curves, and moduli spaces 
$\overline{M}_{0,T}$. The automorphism group of $\overline{M}_{0,T}$
can be canonically identified with the group of bijections $T\to T$.
\smallskip

Now, in order to study the symmetries of the whole operad generated
by $\overline{M}_{0,T}$ we must connect these finite
permutation groups $Aut\,T$ by  families of chosen morphisms
with respect to which on could pass to some meaningful limits.
This is what we did in Section 2, explaining the structure and action
of the permutation group $\S_{\infty}.$

\smallskip

This last and main Section performs this job producing
another ``infinite permutation group'' $\bold{mGT}$ which is a combinatorial
version of the (profinite) Grothendieck--Teichm\"uller group.
We start with showing how to include groupoids of finite sets into a different
poset of groupoids.

\smallskip

We  return  here to the (slightly changed) notations of Sec. 2. Consider the poset $\Cal{N}_*$
whose elements are subsets $\bold{n}:= \{1,\dots , n\}$, $n=1,2,\dots \in \N$,
with binary relation $\bold{m} \le \bold{n}$ iff $\bold{m} \subseteq \bold{n}$.
Obviously, it forms a thin category whose morphisms are cf--maps
coinciding with usual embeddings.

\medskip

Now, extend it to a larger category $\Cal{N}^{cf}_*$ having the same set of objects,
but larger set of morphisms:
$\roman{Hom}_{\Cal{N}^{cf}_*}(\bold{m},\bold{n})$ consists of all cf--maps
obtained by precomposition of a permutation of $\bold{m}$, standard  embedding 
$\bold{m}$ into $\bold{n}$ and postcomposition with a permutation of
$\bold{n}$.

\medskip

{\bf 6.1. Proposition.} {\it (a) $\Cal{N}^{cf}_*$ is a poset in groupoids.

\smallskip

(b) This poset (restricted to $n\ge 3$) acts termwise upon the poset of components
of modular operad $\overline{M}_{0,n}$ whose morphisms are generated by
standard embeddings of locally closed strata.}

\medskip

We will call this poset   $\Cal{N}^{cf}_*$  {\it ``the covariant skeleton of symmetries''}
of our operad.

\medskip

We pass now to the definition of the contravariant skeleton.

\smallskip

Denote by $\Cal{N}^*$ the poset with the same elements as 
$\Cal{N}_*$ but with different ordering: if $\bold{p} =\{1,\dots , p\}$
and $\bold{q} =\{1,\dots , q\}$, then in  $\Cal{N}^*$, $\bold{p}\le \bold{q}$
means that $p$ divides $q$.

\smallskip

Clearly, the map $\Cal{N}^*\to \Cal{N}_*$ identical on elements, is a bijection compatible with respective
order relations, but the inverse map is not compatible.

\smallskip

It is convenient to introduce the family of commutative rings of residues $\Z/q\Z$, $q\ge 3,$
related by the family of natural ring homomorphisms 
$$
t_{q,p}:\, \Z/q\Z \to \Z/p\Z, \quad a\, \roman{mod}\, q \mapsto a\, \roman{mod}\, p
$$
for each pair of natural numbers
$p,q$ such that $p$ divides $q$. It is clear that
if  $p$ divides $q$,$r$ which in turn divide $s$, then
$$
t_{q,p}\circ t_{s,q} = t_{r,p}\circ  t_{s,r} = t_{s,p}.
$$
Now, residue classes of all $d\,\roman{mod}\, q$ with $\roman{g.c.d.}\, (d,q)=1$, 
form the multiplicative group $(\Z/q\Z)^*$. Hence multiplications of $\Z/q\Z$ by them
are permutations that
also act compatibly with all $t_{q,p}$.

\smallskip

We pass now to the central definition of this Section.

\medskip

{\bf 6.2. Definition.} {\it The group $\bold{mGT}_q$ is defined as the subgroup of
permutations of $\Z/q\Z$, generated by the following maps:

\smallskip

(i) multiplications by all elements $d\in (\Z/q\Z)^*$;

\smallskip

(ii) the involution $\theta_q:\, a\mapsto 1-a$.}
\medskip

Notice that $\theta_q$ does not coincide with multiplication
by any $d=d_0$ as above: $\theta_q(0)=1$ whereas $d_0\cdot 0=0$.
\smallskip

Moreover, $\bold{mGT}_q$ is not commutative: $1-da\neq d(1-a)$ in
$\Z/q\Z$ if $d\neq 1$.

\medskip

{\bf 6.3. Proposition.} {\it For each $p, q$ with $p/q$, define the homomorphism
$u_{q,p} :\, \bold{mGT}_q \to \bold{mGT}_p$ by the following
prescription: each permutation of $\bold{q} \in \Cal{N}^*$ belonging to
$\bold{mGT}_q$ is compatible with each map $t_{q,p}$ and after
applying $t_{q,p}$ determines a group homomorphism 
$$
u_{q,p}:\, \bold{mGT}_q\to \bold{mGT}_p
$$
These homomorphism satisfy the following relations:
if  $p$ divides $q$,$r$ which in turn divide $s$, then
$$
u_{q,p}\circ u_{s,q} = u_{r,p}\circ  u_{s,r} = u_{s,p}.
$$
}
This follows directly from the definitions.

\medskip

{\bf 6.4. Corollary.} {\it There exists a well defined group $\bold{mGT}$,
``modified profinite Grothendieck--Teichm\"uller group", which is
the projective limit of groups  $\bold{mGT}_q$ with respect to the
homomorphisms $u_{q,p}$.}

\medskip

It might be still possible to compare $\bold{mGT}$ as an abstract
group given by generators and relations,
with the standard profinite GT--group ([I]) even bypassing
a connection with the Galois group of $\overline{\Q}$
(embedded in $\bold{C}$). 

\smallskip

But actually, such a connection is partially encoded already in our formalism.
Namely, consider the field generated by roots of unity
(i.~e.~the maximal abelian subextension of $\overline{\Q}$ embedded into $\bold{C}$).
Then we can replace each group $\Z/q\Z$ by the group $\mu_q$ of roots
of unity of degree $q$: $a\,\roman{mod}\, q\,\mapsto \, e^{2\pi i a}.$
The action of $(\Z/q\Z)^*$ then becomes the action of the
respective Galois group $e^{2\pi i a}\mapsto e^{2\pi i da}$.
Finally, $\theta_q$  encodes the reflection with respect to $0$ or $\infty$, rather than 1
in  [I].

\smallskip

On the operadic level, the role of indexing families $\mu_q$
in place of $\Z/q\Z$ also becomes more transparent:
imagine $\overline{M}_{0,\mu_q}$ as moduli space 
of deformations of compactified $\bold{G}_m$ with roots of unity, 0, and $\infty$
as marked points. And respective two reflections are related  via conjugation with a very particular element
of  $PSL(2,\Z)$, so the difference is not critical one.

\medskip

{\bf Acknowledgements.} We are very grateful to the Max Planck Institute for Mathematics, for hospitality,
support, and excellent working conditions.

L. Schneps and P. Lochak kindly answered our questions regarding action of the
Grothendieck--Teichm\"uller group, and P. Lochak suggested changes in the 
first draft of this paper ([Lo]).

\medskip

\centerline{\bf References}

\medskip

[BoMa] D. Borisov, Yu. Manin. {\it Generalized operads and their inner
cohomomorhisms.} In: Progress in Math., vol. 265 (2007), pp. 247--308.

\smallskip

[CoMa] N. C. Combe, Yu. I. Manin. {\it Genus zero modular operad and the
profinite Grothendieck--Teichm\"uller group.} Submitted to the
collection dedicated to the memory of Michael Atiyah.

\smallskip

[BrHoRo] P. B. de Brito, G. Horel, M. Robertson. {\it Operads of genus zero curves and the 
Grothendieck--Teichm\"uller group.} Geom. and Topology, 23 (2019), pp. 299--346.

\smallskip

[BrMe] A. Bruno, M. Mella. {\it The automorphism group of $\overline{M}_{0,n}$}.
Journ. Eur. Math. Soc., 15 (2013), pp. 949--968.

\smallskip

[BuSSc] P. Bubenik, V. de Silva, J. Scott. {\it Metrics for generalized persistence modules.}
Found. of Computational Math., 15 (2015), pp. 1501--1531.

\smallskip

[CiMoe] D.-Ch.~Cisinski and I. Moerdijk. {\it Dendroidal Segal spaces and $\infty$--operads.}

Journ. of Topology, 6 (2013), pp. 675--704.
\smallskip

[Fr] B. Fresse. {\it Little discs operads, graph complexes and Grothendieck--Teichm\"uller
groups.} arXiv:1811.12536 .

\smallskip

[HaLoSch] A. Hatcher, P. Lochak, L. Schneps. {\it On the Teichm\"uller tower
of mapping class groups.} J. reine und angew. Math., 521 (2000), pp. 1--24.

\smallskip

[Ho] G. Horel. {\it Profinite completion of operads and the Grothendieck--Teichm\"uller group.}
Advances in Math., 321 (2017), pp. 326--390.

\smallskip

[I] Y.~Ihara. {\it On the embedding of $Gal(\overline{\Q}/\Q)$ into $\widehat{GT}$}.
In [SchGrT], pp. 289--305.

\smallskip

[Ka]  M. Kapranov. {\it Veronese curves and Grothendieck--Knudsen moduli spaces
$\overline{M}_{0,n}$.} Journ. Alg. Geom., 2 (1993), pp. 239--262.

\smallskip

[KashSch] M. Kashiwara, P. Schapira. {\it Categories and Sheaves.}
Grundlehren der math. Wissenschaften, vol. 332, Springer Verlag 2006, 
X + 477 pp.

\smallskip

[Ke] S. Keel. {\it Intersection theory of moduli spaces of stable $n$--pointed curves
of genus zero.} Trans. AMS, 30 (1992), pp. 545--584.

\smallskip

[KoMa] M. Kontsevich, Yu. Manin. {\it Gromov-Witten classes, quantum cohomology, and enumerative geometry.}
Comm. Math. Phys., 164 (1994),  pp. 525--562.

\smallskip

[Lo] P. Lochak. {\it Letters of 20.05.2019 and 12.05.2019}.

\smallskip

[MaMar] Yu. Manin, M. Marcolli. {\it Nori diagrams and persistent homology.}
Submitted to Mathematics in Computer Science. arXiv:1901.1031. 45 pp

\smallskip

[SchGrT] L.~Schneps, ed. {\it The Grothedieck Theory of Dessins d'enfants.}
London Math. Soc. Lecture Note Series 200. Cambrifge UP, 1994.

\smallskip

[Sch] L.~Schneps. {\it Dessins d'enfants on the Riemann sphere.}
In [SchGrT], pp. 47--77.

\bigskip
{\bf No\'emie C. Combe,  Max Planck Institute for Mathematics, Vivatsgasse 7, 53111
Bonn

\medskip
Yuri I. Manin, Max Planck Institute for Mathematics, Vivatsgasse 7, 53111
Bonn}

\enddocument